\documentclass[12pt]{amsart}

\usepackage{amsmath,amsthm,amssymb}
\usepackage[T1]{fontenc}
\usepackage{latexsym}
\usepackage{subfigure, epsfig, epsf}
\usepackage{color,fancyhdr,lastpage,graphicx}
\usepackage{xspace}
\usepackage{wasysym}
\usepackage{setspace}
\usepackage{bbm}

\newcommand{\HH}{\mathbb{H}}        \newcommand{\RR}{\mathbb{R}}          
                \newcommand{\MMM}{\mathbb{M}}
\newcommand{\OO}{\mathbb{O}}        \newcommand{\VV}{\mathbb{V}}

 \newcommand{\mcC}{\mathcal{C}} \newcommand{\mcF}{\mathcal{F}} \newcommand{\mcU}{\mathcal{U}}
 \newcommand{\mcV}{\mathcal{V}}  
    \newcommand{\mcG}{\mathcal{G}}
 \newcommand{\mcE}{\mathcal{E}}

  \newcommand{\be}{{\bf e}}  
    
 \newcommand{\bff}{{\bf f}}

 \newcommand{\wbfT}{\widetilde{{\bf T}}}

\newcommand{\ra}{\rightarrow}

\newcommand{\EE}{\mathbb{E}} \newcommand{\PP}{\mathbb{P}}

\newcommand{\bm}{Brownian motion }

\newcommand{\rv}{random variable }

\newcommand{\sde}{stochastic differential equation }

\newcommand{\as}{almost-surely }

\newcommand{\Vol}{\small{\textsc{Vol}}}

\newcommand{\Ric}{\textsc{Ric}}  \newcommand{\wRic}{\widetilde{\textsc{Ric}}}

\newcommand{\st}{such that }

\newcommand{\ep}{\epsilon}

\renewcommand{\leq}{\leqslant}
\renewcommand{\geq}{\geqslant}

\newcommand{\la}{\lambda}
\newcommand{\al}{\alpha}

\newcommand{\wrt}{with respect to }
\renewcommand{\st}{such that }

\newcommand{\ssk}{\smallskip}

\newcommand{\noi}{\noindent}

\newtheorem{thm}{ \hspace{-0.15cm}{\sc Theorem} }
\newtheorem{lem}[thm]{ \hspace{-0.15cm}{\sc Lemma} }
\newtheorem{prop}[thm]{ \hspace{-0.15cm}{\sc Proposition}}

\numberwithin{equation}{section} 

\newenvironment{Dem}{%
    \begin{list}{$\quad \,$ {\sc Proof --}}{%
        \setlength{\topsep}{0pt}%
        \setlength{\leftmargin}{0pt}%
        \setlength{\rightmargin}{0pt}%
        \setlength{\listparindent}{0pt}%
        \setlength{\itemindent}{0pt}%
        \setlength{\parsep}{0pt}%
        \addtolength{\leftmargin}{20pt}%
        \addtolength{\rightmargin}{0pt}%
    } \item }{\hfill{\space $\rhd$}\end{list}\smallskip}

\title[A probabilistic view on singularities]{A probabilistic view on singularities}
\date{\today}
\author[I. Bailleul]{Isma\"el Bailleul}
\address{Statistical Laboratory, Center for Mathematical Sciences, Wilberforce Road, Cambridge, CB3 0WB, UK}
\email{i.bailleul@statslab.cam.ac.uk}
\urladdr{http://www.statslab.cam.ac.uk/~ismael/}


\begin{document}

\begin{abstract}
The aim of this article is to promote the use of probabilistic methods in the study of problems in mathematical general relativity. Two new and simple singularity theorems, whose features are different from the classical singularity theorems, are proved using probabilistic methods. Under some energy conditions, and without any causal or initial/boundary assumption, simple conditions on the energy flow imply probabilistic incompleteness. 
\end{abstract}

\maketitle

\begin{center}
\tableofcontents
\end{center}

\section{Introduction}

The introduction of the concepts of singularities and boundary of a spacetime both have their source in the observation that Einstein's geometrical picture of a spacetime does not prevent the existence of some undesirable features like the existence of regions where some geometrical scalar explodes along some path, or the existence of physically relevant incomplete paths (geodesics, paths with bounded acceleration, etc). Following the pioneering works of Penrose \cite{Penrose} and Hawking \cite{Hawking2}, \cite{Hawking1}, \cite{HawkingPenrose}, most singularity theorems state that a spacetime has an incomplete causal geodesic  provided some energy, causality and boundary or initial conditions hold. We consider in this article a purely geometrical random dynamics producing random timelike paths, to be considered as randomly perturbated geodesics. Technically speaking, we are going to construct some non-trivial probabilities on the (separable metric) space of inextendible timelike paths parametrized by their proper time, with finite or infinite lifetime. We prove in section \ref{SectionMainThm} that some energy conditions and conditions on the energy flow are sufficient to ensure the existence of a set of inextendible incomplete timelike paths of positive probability. This provides an unusual type of conclusion under a non-common set of conditions where no causality assumption is required, nor any boundary or initial condition. 

\medskip

Cartan's moving frame language provides a very intuitive way of describing $\mcC^2$ trajectories by transporting parallelly an initial frame along the path and describing the variations of the velocity in that moving frame. The datum of these vector space-valued variations suffices to reconstruct the path. Taking them as the primary object and choosing them randomly produces random paths (forgetting all possible technical problems). Our construction of random timelike paths relies on a variant of this procedure and associates to any starting point $(m_0,\dot m_0)$ in the future-oriented unit bundle $T^1\MMM$ a timelike path in $\MMM$ under the form of a continuous path $(m_s,\dot m_s)$  in $T^1\MMM$ subject to the condition $\frac{d}{ds}m_s = \dot m_s$. So our random dynamics is a dynamics in $T^1\MMM$, in the same way as the (timelike) geodesic motion is better viewed as a flow in the unit bundle. The (starting point dependent) distribution of these random timelike paths will provide the above mentionned set of probabilities.

\medskip

These random dynamics are called ``diffusions'' in the sequel. The reader should have in mind that these diffusionss are not physical diffusions in the sense that they are not models for diffusion phenomena in a relativistic medium like a gaz. We refer for this kind of questions to the works by Debbasch and his co-workers, \cite{ROUP0}, \cite{ROUPCurved}, or Dunkel and Hanggi, \cite{DunkelHanggi}, \cite{DunkelHanggiRW}, their review \cite{DunkelHanggiReview}, and the references cited therein. Rather, the diffusions considered below are what probabilists call a diffusion: intuitively, a random process $X = (X^1,\dots,X^n)$ \st
\begin{equation*}
\begin{split}
&\EE\bigl[X^i_{s+\ep}-X^i_s\bigr] = b^i(X_s)\ep + o(\ep), \\
&\EE\bigl[(X^i_{s+\ep}-X^i_s)(X^j_{s+\ep}-X^j_s)\bigr] = a^{ij}(X_s)\ep + o(\ep).
\end{split}
\end{equation*}
for $i,j=1..n, s\geq 0$ and $\ep>0$. Here, $\mcF_s = \sigma(X_r\,;\,r\leq s)$ and $b^i$ and $a^{ij}$ are measurable functions; see chapter 5 in \cite{RW}.

\medskip

Section \ref{SectionRelDiff} gives a brief overview of the diffusions considered here, with the help of which we derive two simple probabilistic singularity theorems in section \ref{SectionMainThm}, theorems \ref{MainTheorem} and \ref{2ndMainTheorem}. As mentionned above, the random dynamics considered in this work are random perturbations of the geodesic flow in $T^1\MMM$. The main issue addressed in this work is the following question: \textit{Can (timelike) geodesic completeness be destroyed by a random perturbation of its dynamics?} Theorems \ref{MainTheorem} and \ref{2ndMainTheorem} both provide conditions under which the answer to this question is positive. The core of their proofs is a stochastic analogue of the following kind of trivial observation about the unperturbed flow. If there exists two functions $f\leq h$ and an initial condition $\phi_0\in T^1\MMM$ \st $f$ is bounded below by $e^{c\,s}$ and $h$ bounded above by $e^{c's}$ along the geodesic started from $\phi_0$, then this geodesic cannot be future complete if $c'<c$. The stochastic counterpart of this observation is subtler though.

\medskip

As this work is the first ever written on the subject, we have chosen to present only some basic aspects of the situation. We hope this will convince the reader of the possible intest of using probabilistic methods in problems about spacetime geometry.

\section{Relativistic diffusions}
\label{SectionRelDiff}

\noi One owes to Dudley \cite{Dudley1} and Schay \cite{Schay1}, \cite{Schay2}, the merit to have asked first if there is a natural way of defining random timelike paths in Minkowski spacetime, without introducing any additional structure than the metric. The models of physical diffusions considered e.g. by Debbasch et al. \cite{ROUP0}, Angst and Franchi \cite{JurgenFranchi}, or Dunkel and Hanggi \cite{DunkelHanggi}, require an additional vector field to be defined. Dudley provided a complete answer to the above problem but could not pursue further his investigations to the general relativistic framework for lack of technical tools available at that time. The subject remained untouched for nearly forty years before Franchi and Le Jan generalized Dudley's motion to a general Lorentzian manifold in \cite{FranchiLeJan}. We give a brief overview of these dynamics in section \ref{SectionRelDiffMin} and \ref{SectionRelDiffLorentz} below as they will be our main tool to probe the existence of spacetime probabilistic singularities in section \ref{SectionMainThm}.

\subsection{Relativistic diffusion in Minkowski spacetime}
\label{SectionRelDiffMin}

The question asked by Dudley and Schay is the following: Is there a natural way of defining Markovian random timelike paths, independently of any reference frame, and without using any other object than Minkowski metric? To describe their answer, denote by $g$ Minkowski metric and write $\HH = \bigl\{\dot{m}\in\RR^{1,3}\,;\,\dot{m}^0>0,\,g(\dot{m},\dot{m})=1\bigr\}$ for the upper half-sphere. Although $g$ has signature $(-,+,+,+)$, its restriction to any tangent space of the spacelike hypersurface $\HH$ is definite-positive, so $\HH$ inherits from the ambient space a Riemannian structure which actually turns it into a model of the $3$-dimensional hyperbolic space. \bm on $\HH$ is defined as the unique continuous Markov process $\{\dot{m}_s\}_{s\geq 0}$ with generator half of the Laplacian $\triangle$ of $\HH$. (We have indexed \bm by $\RR_+$; it can indeed be proved that it has an infinite lifetime; see e.g. \cite{Pinsky}.) As any timelike path $\{m_s\}_{s\geq 0}$ indexed by its proper time is differentiable at almost-all times, and determined by its derivative $\dot{m}_s$, since $m_s = m_0 + \int_0^s\dot{m}_u\,du$, a random timelike path is determined by an $\HH$-valued random process $\{\dot{m}_s\}_{s\geq 0}$. So, if one wants to talk of a random timelike paths-valued Markov process, one needs to record the position and the velocity of the path in the state space and work in $\RR^{1,3}\times\HH$. Dudley showed in \cite{Dudley1} that there exists essentially a unique $\RR^{1,3}\times\HH$-valued continuous Markov process $\bigl\{(m_s,\dot{m}_s)\bigr\}_{s\geq 0}$ \st 
\begin{itemize}
   \item one has $m_s = m_0 + \int_0^s\dot{m}_u\,du$, for all $s\geq 0$,
   \item its law is invariant by the action of the affine isometries of $\RR^{1,3}$.
\end{itemize}
It corresponds to a Brownian velocity process $\{\dot m_s\}_{s\geq 0}$ with generator $\frac{\sigma^2}{2}\,\triangle$, where $\sigma$ is a positive constant. In short, there is (essentially) a unique way of constructing a random timelike path, by imposing to the velocity to undergo Brownian oscillations in $\HH$.

\ssk

It will clarify the construction in the general framework of a Lorentzian manifold to give a slightly different picture of Dudley's process. Consider the random motion in spacetime of an infinitesimal rigid object, represented by a path in $\RR^{1,3}\times SO_0(1,3)$, where we write $SO_0(1,3)$ for the identity's component in $SO(1,3)$. Denoting by $(m,e)=\bigl(m,(e_0,e_1,e_2,e_3)\bigr)$ a generic element of $\RR^{1,3}\times SO_0(1,3)$, the map $\pi_1\bigl((m,e)\bigr) = (m,e_0)$ is a projection from $\RR^{1,3}\times SO_0(1,3)$ onto $\RR^{1,3}\times\HH$. Denote by $\{\ep_0,\ep_1,\ep_2,\ep_3\}$ the canonical basis of $\RR^{1,3}$. For each $j\in\{1,2,3\}$, the Lie element $E_j = \ep_0\otimes \ep_j^* + \ep_j\otimes \ep_0^*\in so(1,3)$ generates a hyperbolic rotation in the $2$-plane spanned by $\ep_0$ and $\ep_j$. Define on $SO_0(1,3)$ some left invariant vector fields $V_j$ setting
$$
V_j(e) = eE_j, \quad j\in\{1,2,3\}.
$$
Let $w$ be a $3$-dimensional Brownian motion and use the notation ${\circ}d$ for Stratonovich differential. Consult section 5, chap. V, of the book \cite{RW}, or chapter V of \cite{IW}, for the necessity of using Stratonovich formalism when dealing with stochastic differential equations on manifolds, and for a comparison with Ito's differential. By construction, the solution $(e_s)_{s\geq 0}$ of the Stratonovich stochastic differential equation on $SO_0(1,3)$
\begin{equation}
\label{SdeBmSO}
{\circ d}e_s = \sum_{j=1}^3V_j(e_s)\,{\circ d}w^j_s
\end{equation}
projects down by $\pi_1$ into a \bm on $\HH$. (Consult the book \cite{ElworthyBook} of Elworthy, or \cite{IW}, if you do not feel comfortable with this fact.) So Dudley's dynamics is described by equation \eqref{SdeBmSO} and
\begin{equation}
\label{SdeDudleySpaceTime}
dm_s = e_0(s)\,ds.
\end{equation}

\ssk

\noi \textbf{Remark.} It is interesting to notice that Dudley's diffusion is obtained heuristically as the large scale limit of geodesics in the discrete random approximations of a spacetime introduced by Sorkin in his causal set theory; see \cite{DowkerHensonSorkin}, \cite{PhilpottDowkerSorkin} and \cite{PhilpottSimulation}.

\subsection{Basic relativistic diffusion on a general Lorentzian manifold}
\label{SectionRelDiffLorentz}

Let $(\MMM,g)$ be a Lorentzian manifold, oriented and time-oriented. (These assumptions are harmless as they hold on a finite covering of $\MMM$.) Denote by $T^1\MMM$ the unit future-oriented bundle over $\MMM$ and by $\OO\MMM$ the component of the orthonormal frame bundle made up of pairs $(m,\be)$, with $m\in\MMM$ and $\be=(\be_0,\be_1,\be_2,\be_3)$ a direct orthonormal basis of $T_m\MMM$ with $\be_0$ future-oriented. We shall denote by $\phi = (m,\dot m)$ a generic point of $T^1\MMM$ and by $\Phi = (m,\be)$ a generic element of $\OO\MMM$. Write $\pi_1 : \OO\MMM\ra T^1\MMM$ and $\pi_0 : \OO\MMM\ra \MMM$ for the canonical projections.

\medskip

\noi As can be expected from the Minkowskian picture, the basic relativistic diffusion on a Lorentzian manifold $\MMM$ is actually a diffusion in $T^1\MMM$. Roughly speaking, one can construct this process by rolling without splitting the trajectories of the relativistic diffusion in Minkowski spacetime on $\MMM$. A more formal way of proceeding is to introduce an $\OO\MMM$-valued $SO(3)$-invariant diffusion process whose projection in $T^1\MMM$ is consequently a diffusion process on its own.

\ssk

The choice of $\OO\MMM$ as a framework is motivated by the fact that it bears more structure than $T^1\MMM$ and is the natural framework where to use Cartan's ideas on moving frames. The action of $SO_0(1,3)$ on each fiber of $\pi_0$ induces the canonical vertical vector fields; denote by $V_j$ the vector field associated with the Lie element $E_j$. Denote by $(H_i)_{i=0..3}$ the canonical horizontal vector fields on $\OO\MMM$ associated with Levi-Civita connection. The \textbf{basic relativistic diffusion} is defined in a dynamical way as the unique solution of the following Stratonovich \sde on $\OO\MMM$
\begin{equation}
\label{SdeRelDiff}
{\circ d}\Phi_s = H_0\bigl(\Phi_s\bigr)\,ds + \sigma\sum_{j=1}^3 V_j\bigl(\Phi_s\bigr)\,{\circ d}w^j_s,
\end{equation}
where $w$ is a $3$-dimensional \bm and $\sigma$ a positive constant. It describes a random perturbation of the geodesic flow whose intuitive meaning is the following. To get $\Phi_{s+ds}$ out of $\Phi_s=(m_s,\be_s)$, transport first $\be_s$ parallelly along the geodesic starting from $m_s$ in the direction $\be_0(s)$, during an amount of time $ds$; you get an orthonormal frame $\bff_s$ of $T_{m_{s+ds}}\MMM$. Making then, in each spacelike $2$-plane of $\bff_s$, independent hyperbolic rotations of angle a centered normal \rv with variance $\sigma^2\,ds$, you get $\be_{s+ds}$. This dynamics is the straightforward generalization of Dudley's dynamics \eqref{SdeBmSO} and \eqref{SdeDudleySpaceTime} and was first considered in \cite{FranchiLeJan}. The following statement gives a different view on this dynamics.

\begin{lem}
\label{LemmaLifting}
Let $\gamma : [0,T]\ra\MMM$ be a $\mcC^2$ timelike path parametrized by its proper time, and $\Gamma_0\in\OO\MMM$ \st $\pi_1(\Gamma_0) = \bigl(\gamma(0),\dot\gamma(0)\bigr)\in T^1\MMM$. Then there exists a unique $\mcC^2$ path $\bigl(\Psi_s\bigr)_{0\leq s\leq T}$ in $\OO\MMM$, and some unique $\mcC^1$ real-valued controls $h^1,h^2,h^3$ defined on $[0,T]$, \st $\Psi_0 = \Gamma_0,\;\pi_1(\Psi_s) = \bigl(\gamma(s),\dot \gamma(s)\bigr)$ and 
$$
\dot \Psi_s = H_0(\Psi_s) + \sum_{j=1}^3\limits V_j(\Psi_s)\,h^j(s).
$$
\end{lem}
\noi So the relativistic diffusion dynamics is obtained by replacing the deterministic controls of a typical $\mcC^2$ timelike path by Brownian controls.

\ssk

Equation \eqref{SdeRelDiff} has a unique strong solution defined up to its explosion time $\zeta$. Write $\OO\MMM$ as an increasing union of relatively compact open sets $O_n$ and define the stopping times $T_n=\inf\{s\geq 0\,;\,\Phi_s\notin O_n\}$. The explosion time $\zeta$ is by definition  the increasing limit of the $T_n$, and does not depend on the arbitrary choice of sets $O_n$. It is not difficult to see on this equation that it defines an $SO_0(3)$-invariant diffusion which, as a consequence, projects on $T^1\MMM$ into a diffusion process. Consult \cite{FranchiLeJan}, theorem 1, or \cite{BailleulIHP}, theorem 4, for the details. The diffusion $\bigl(\Phi_s\bigr)_{0\leq s<\zeta}$ in $\OO\MMM$ has generator
$$
\mcG = H_0 + \frac{\sigma^2}{2}\sum_{j=1}^3 V_j^2,
$$
where we consider vector fields as first order differential operators. In so far as we are primarily interested in the $T^1\MMM$-valued process (as it provides us directly with random timelike paths), while we shall mainly be working with the above $\OO\MMM$-valued process it is important to notice that 

\begin{prop}[\cite{BailleulFranchi}, prop. 1]
The $\OO\MMM$-valued diffusion $\bigl(\Phi_s\bigr)_{0\leq s<\zeta}$ and its $\pi_1$ projection in $T^1\MMM$ have the same lifetime.
\end{prop}

We shall thus freely work in the sequel with the $\OO\MMM$-valued diffusion $(\Phi_s)_{0\leq s<\zeta}$.

\medskip

\noi \textbf{Remarks. (i)} You might ask why we called \textit{basic} relativistic diffusion the solution of equation \eqref{SdeRelDiff}, and not simply relativistic diffusion. This is due to the fact that, contrary to what happens in Minkowski spacetime, where there is only one way of constructing nice random timelike paths, there are many ways of doing it on a general Lorentzian manifold. Think for example of a diffusivity $\sigma$
in \eqref{SdeRelDiff} depending on the location of the particle (it may be the scalar curvature of the manifold at that point for instance). Consult \cite{BailleulIHP} and \cite{FranchiLeJanCurvature} for more material on this subject, seen from a mathematical point of view, and \cite{ROUPUnifying2} for a physical point of view on related matters. Let us repeat here that this diffusion is not to be thought of as a mathematical model for a physical diffusion phenomenon but rather as a mathematical object useful for studying some features of the spacetime geometry.

\ssk

\textbf{(ii)} Basic relativistic diffusions have only been studied explicitl in a few cases: in Minkowski spacetime \cite{BailleulPoisson}, Robertson-Walker spacetimes \cite{Angst}, Schwarzschild \cite{FranchiLeJan} and G\"odel \cite{Franchi} spacetimes. The stochastic completeness question is easily delt in each case. While it is trivial to find a geodesically incomplete spacetime which is stochastically complete (remove a point from Minkowski spacetime -- see paragraph \ref{SectionIncompleteness} below), the possibility to have a (timelike) geodesically complete and stochastically incomplete spacetime has noot been established so far; theorems \ref{MainTheorem} and \ref{2ndMainTheorem} address that issue.

\ssk

\textbf{(iii)} As mentionned in the introduction, one can consider equation \eqref{SdeRelDiff} as a dynamical way of constructing a probability measure on the separable metric space of inextendible timelike paths parametrized by their proper time (with finite or infinite lifetime). Given a starting point $\Phi_0\in\OO\MMM$, the distribution $\PP_{\Phi_0}$ of $(\Phi_s)_{0\leq s<\zeta}$ is such a probability. Each such probability is highly non-trivial, as shown by the following qualitative statement. Let $\gamma : [0,T] \ra \MMM$ be a $\mcC^2$ timelike path parametrized by its proper time, and $\bigl(\be_1(0),\be_2(0),\be_3(0)\bigr)$ be \st $\bigl(\dot \gamma(0),\be_1(0),\be_2(0),\be_3(0)\bigr)$ is a frame of $T_{\gamma(0)}\MMM$. Transporting parallelly $\bigl(\be_1(0),\be_2(0),\be_3(0)\bigr)$ along $\gamma$ provides a frame $\bigl(\dot \gamma(s),\be_1(s),\be_2(s),\be_3(s)\bigr)$ of $T_{\gamma(s)}\MMM$.  The map 
\vspace{-0.15cm}$$
F : (s,x)\in [0,T]\times\RR^3\mapsto\exp_{\gamma(s)}\Bigl(\sum_{j=1}^3x^j\be_j(s)\Bigr)
\vspace{-0.15cm}$$ 
is a well-defined diffeomorphism from $[0,T]\times U$ onto its image $\mcV$, for some small enough open ball $U$ of $\RR^3$ with center $0$; $\mcV$ is a tube around $\gamma$.

\begin{prop}
\label{PropTube}
Let $\Phi=(m,\be)\in\OO\MMM$ be \st $\gamma(0)\in I^+(m)$. Then the $\PP_{\Phi}$-probability that the basic relativistic diffusion hits $F\bigl(\{0\}\times U\bigr)$ and exits the tube $\mcV$ in $F\bigl(\{T\}\times~ U\bigr)$ is positive.
\end{prop}

\medskip

\noi \textbf{An example: basic relativistic diffusion in Schwarzschild spacetime.} Let denote by $(\MMM,g)$ Kruskal-Szekeres extension of Schwarzschlid spacetime. Franchi and Le Jan proved in their seminal work \cite{FranchiLeJan} that, for any starting point, the relativistic diffusion hits the boundary of the black hole with a positive probability strictly smaller than $1$. After that time, as any other timelike path, the trajectory of the diffusion hits the singularity before $\frac{\pi}{2}R$ units of proper time have ellapsed, where $R$ is the radius of the black hole. So $\zeta$ is finite with positive probability, whatever the starting point of the diffusion.

\section{Probabilistic incompleteness}
\label{SectionMainThm}

\subsection{Probabilistic incompleteness and geometry}
\label{SectionIncompleteness}

A probabilistic incompleteness theorem is a statement of the form $\PP_{\Phi}(\zeta<\infty)>0$, for some initial starting point $\Phi\in\OO\MMM$ of the diffusion process. It should be clear from the description of Dudley's diffusion given in section \ref{SectionRelDiffMin} that probabilistic incompleteness and geodesic incompleteness are two different notions in general. Indeed, removing a point from Minkowski spacetime produces a geodesically incomplete spacetime while Dudley's diffusion is \as defined for all proper times. This is easily seen as follows. Choose a frame and denote by $t_s$ the time component of the $\RR^{1,3}$-part of the diffusion. As $\frac{d}{ds}t_s\geq 1$, one can change parameter and use the time of the frame rather than the proper time to parametrize the trajectory. This amounts to consider the diffusion as an $\bigl(\RR^3\times\HH\bigr)$-valued diffusion. It is easily seen to be hypoelliptic, using H\"ormander's criterion, so that its distribution at any given time has a density. If the point removed belongs to the slice $\{t=T\}$, the diffusion has a null probability of hitting it. It is also true that any $2$-dimensional submanifold of $\OO\MMM$ is polar for the basic relativistic diffusions. We describe in theorems \ref{MainTheorem} and \ref{2ndMainTheorem} below geodesically complete spaces which are probabilistically incomplete.

\ssk

In a Riemannian setting, it is well-understood that any statement about the Laplacian, its eigenvalues/eigenfunctions, heat kernel, zeta function, etc. provides a non-elementary information about the geometry of the manifold. In so far as the generator $\mcG$ of the basic relativistic diffusion is constructed out of canonical geometric vector fields, any statement about $\mcG$ provides a non-elementary information about the geometry of $\OO\MMM$, or $T^1\MMM$.

The following proposition links probabilistic incompleteness and geometry; its proof is similar to the proof of the corresponding statement for the Laplacian, and can be found e.g. in theorem 6.2. of \cite{Grigoryan}.

\begin{prop}
\label{PropEquivalencePDE}
The following statements are equivalent.
\begin{enumerate}
   \item There exists a point $\Phi\in\OO\MMM$ \st $\PP_{\Phi}(\zeta<\infty)>0$.
   \item Let $\la>0$. There exists a non-null bounded function $f : \OO\MMM\ra\RR$ \st \\ $(\mcG-\la)f = 0$.
   \item Let $T>0$. There exists a non-null bounded solution to the Cauchy problem \\ $\partial_t h = \mcG h$, on $[0,T]\times\OO\MMM$, with initial condition $0$.
\end{enumerate}
\end{prop}

\noi This equivalence between the probabilistic incompleteness problem and the above two problems on \textit{linear} partial differential equations brings a different point of view on the probabilistic explosion problem, and a huge tool kit to investigate it. There is no similar correspondence for the classical geodesic incompleteness problem. The following proposition clarifies what happens in case of explosion.

\begin{prop}
Suppose there exists a $\Phi\in\OO\MMM$ \st $\PP_\Phi(\zeta<\infty)>0$.Then explosion occurs in an arbitrarily small time with positive probability: $\PP_\Phi(\zeta<\ep)>0$, for all $\ep>0$.
\end{prop}

\begin{Dem}
Set $f(\Psi,s) = \PP_\Psi(\zeta<s)$, for $\Psi\in\OO\MMM$ and $s>0$; due to H\"ormander's theorem, this function depends smoothly on $(\Psi,s)\in\OO\MMM\times\RR_+^*$. Write $\Phi=(m,\be)$ and set $I^+(\Phi) = \pi^{-1}\bigl(I^+(m)\bigr)$. Given $\ep>0$, suppose $\PP_\Psi(\zeta<\ep)=0$ for all $\Psi\in I^+(\Phi)$. Using the strong Markov property first, then the Markov property inductively and the fact that $I^+(\Psi)\subset I^+(\Phi)$, for $\Psi\in I^+(\Phi)$, we would have $\PP_\Phi(\zeta\geq k\ep)=1$, for all $k\geq 1$, contradicting $\PP_\Phi(\zeta<\infty)>0$. So there exists an element $\Psi_0\in I^+(\Phi)$ \st $f(\Psi_0,\ep)>0$\,; by continuity, this remains true in an open neighbourhood of $\Psi_0$. As, by proposition \ref{PropTube}, one can reach this open set in an arbitrarily small time from $\Phi$ with positive probability, the result follows from the strong Markov property.
\end{Dem}

\ssk

Probabilistic incompleteness results have two other noticeable distinct features when compared with the classical deterministic singularity theorems. As a function of $\Phi$, the quantity $\PP_{\Phi}(\zeta<\infty)$ is known to be $\mcG$-harmonic ($\mcG f = 0$) and satisfies the strong maximum principle. So if $\PP_{\Phi}(\zeta<\infty)=1$ at some point $\Phi=(m,\be)\in\OO\MMM$, we actually have $\PP_{\Phi'}(\zeta<\infty)=1$ for any $\Phi'=(m',\be')$ with $m'$ in the chronological future of $m$.

\ssk

Given a spacelike hypersurface $\VV$ of $\MMM$ denote by $\sigma_{\VV}$ the volume measure on $\VV$ inherited from the ambient geometry and set $\OO\VV=\bigl\{\Psi=(m,\be)\in\OO\MMM\,;\,m\in\VV\bigr\}$. Write $\Vol_{\OO\VV}$ for the natural volume measure $\Vol_m(d\be)\otimes\sigma_\VV(dm)$ on $\OO\VV$, where $\Vol_m(d\be)$ is the Haar measure on the fiber $SO_0(1,3)$ above $m$. Given a point $\Phi_0=(m_0,\be_0)\in\OO\MMM$, there exists a relatively compact neighbourhood $\mcU$ of $m_0$ in $\MMM$ \st $\mcU$ is included in a larger globally hyperbolic neighbourhood of $m_0$ and $\partial\bigl(\mcU\cap I^+(m_0)\bigr)$ is the union of a lightlike hypersurface and a smooth spacelike hypersurface $\VV$. Denote by $H$ the random hitting time of $\OO\VV$ by $(\Phi_s)_{0\leq s<\zeta}$; it is $\PP_{\Phi_0}$-\as finite. The hypoellipticity of the generator $\mcG$ of the diffusion garantees that the distribution of $\Phi_T$ has a smooth positive density $D(\Psi)$ \wrt $\Vol_{\OO\VV}$. Consult proposition 5 of \cite{BailleulIHP} for a proof. We can thus write 
$$
\PP_{\Phi}(\zeta<\infty) = \EE_{\Phi}\bigl[\PP_{\Phi_H}(\zeta<\infty)\bigr] = \int_{\OO\VV}\PP_{\Psi}(\zeta<\infty)\,D(\Psi)\Vol_{\OO\VV}(d\Psi).
$$
So if $0<\PP_{\Phi}(\zeta<\infty)$, there is a subset $\mcE$ of $\OO\VV$ of positive $\Vol_{\OO\VV}$-measure \st $\PP_{\Psi}(\zeta<\infty)>0$ for all $\Psi\in\mcE$. Not only do we have a set of positive probability of incomplete inextendible timelike paths with a common starting point, but this conclusion also holds for a non-trivial uncountable collection of starting points.

\subsection{Completeness conditions}
Despite all the works done, it remains unclear what precise features of a spacetime forbid the existence of honest inextendible incomplete timelike paths. The example found by Geroch \cite{GerochSingularity} of a timelike and lightlike geodesically complete spacetime having an inextendible incomplete path with bounded acceleration gives an idea of the subtleties involved in that matter (see also \cite{BeemSingularHyperbolic}).

\medskip

As far as our method is concerned, the recent work \cite{BailleulFranchi} delimitates its domain of application by determining some general situations in which the paths of the diffusion have \as an infinite lifetime. We recall here two such conditions and refer the reader to that work for more and subtler material.

\ssk 

$\bullet$ Let $\MMM = I\times S$ be a globally hyperbolic spacetimes whose metric tensor is of the form 
$$
g_m(q,q) = -a_m^2\,\big|q^0\big|^2 + h_{m}(q^S,q^S), \quad q\in T_m\MMM,
$$
where $q^0$ is the image of $q$ by the differential of the first projection $I\times S \ra I$ and $q^S$ the image of $q$ by the differential of the second projection $I\times S \ra I$. Write $m=(t,x)\in I \times S$. The function $a$ is a positive $\mcC^1$ function and $h_{m}$ a positive-definite scalar product on $T_xS$, depending in a $\mcC^1$ way on $m$. This class of spacetimes contains as elements all Robertson-Walker spacetimes -- in particular de Sitter and Einstein-de Sitter spacetimes -- and the universal covering of the anti-de Sitter spacetime. 

\begin{prop}[\cite{BailleulFranchi}]
The relativistic diffusion does not explode if $\nabla a$ is everywhere non-spacelike, future-oriented.
\end{prop}

\ssk

$\bullet$ It has been proved in \cite{BailleulFranchi} that if the spacetime is b-complete in the sense of Schmidt (a strong requirement; see e.g. \cite{SchmidtbComplete1} and \cite{SchmidtbComplete2}) then the relativistic diffusion does not explode.

\subsection{A first probabilistic singularity theorem}
\label{SectionFirstTheorem}


\noi \textbf{a) What can we expect?} Pursuing the seminal works of Penrose and Hawking, the proofs of most singularity theorems follow one of the following two lines of reasonning (consult Senovilla's review \cite{Senovilla} for a thorough and critical review on singularity theorems, or Wald's book \cite{Wald}).
\begin{enumerate}
   \item Using causality conditions, one constructs a useful maximal geodesic. At the same time, the energy and initial/boundary conditions induce conjugate or focal points along any geodesic, provided time is allowed to run long enough; as this cannot happen along a maximal geodesic, its time parameter has to be bounded.
   \item Supposing the spacetime geodesically complete, one constructs a compact proper achronal boundary whose existence prevents the existence of an open Cauchy hypersurface.
\end{enumerate}
Each of these schemes uses in a crucial way the rich structure of the geodesic flow. In so far as the random flow generated by equation \eqref{SdeRelDiff} is not easy to analyse, it is not obvious to see which geometric features of a spacetime can lead to exploding solutions of equation \eqref{SdeRelDiff}.

\medskip

The following very general and rough explosion result will be sufficient for our purpose in the next section; it is due to Khasminski. A generalisation of this result is proved in lemma \ref{2ndExplosionLemma}. Recall the equivalence between points $(1)$ and $(2)$ in proposition \ref{PropEquivalencePDE}.

\begin{lem}
\label{LemmaExplosion}
Let $f :\OO\MMM\ra\RR$ be a $\mcC^2$ function bounded above by a positive constant and $\Phi_0$ a point of $\OO\MMM$ \st $f(\Phi_0)>0$. If there exists a positive constant $C$ \st $\mcG f\geq Cf$, then $\PP_{\Phi_0}\bigl(\zeta<\infty\bigr)>0$.
\end{lem}

\medskip

\noi \textbf{b) A probabilistic incompleteness theorem.} Let $R(\cdot,\cdot)\cdot$ be Riemann curvature tensor and denote by $\Ric$ the Ricci curvature tensor and by $\wRic$ its restriction to $T^1\MMM:\;\; \wRic(\phi) := \Ric_m(\dot m,\dot m)$, for $\phi=(m,\dot m)\in T^1\MMM$. We identify $\wRic$ to a function on $\OO\MMM$ setting $\wRic(\Phi) = \wRic\bigl(\pi_1(\Phi)\bigr)$. Write $R$ for the scalar curvature and ${\bf T} = \Ric-\frac{1}{2}\,R\,g$ for the energy-momentum tensor.

\begin{thm}
\label{MainTheorem}
Let $(\MMM,g)$ be a Lorentzian manifold. Suppose 
\begin{enumerate}
   \item $\forall\,(m,\dot m)\in T^1\MMM, \;\; {\bf T}_m(\dot m,\dot m)\geq 0$,
   \item the function $\wRic$ is non-identically constant, positive at some point $\Phi_0$, and bounded above,
   \item there exists a constant $c\in\RR$ \st $\;H_0\,\wRic \geq c\,\wRic$.
\end{enumerate}
Write $c= C-2\,\sigma^2$ for some positive constants $C$ and $\sigma$. Then the basic relativistic diffusion with diffusivity $\sigma^2$ explodes with $\PP_{\Phi_0}$-positive probability. 
\end{thm}

\ssk

The proof makes a crucial use of the following simple lemma.

\begin{lem}
\label{LemmaGenerator}
We have $\;\;\mcG\,\wRic = H_0\,\wRic + 2\,\sigma^2\,\wRic + 2\,\sigma^2\,\widetilde{{\bf T}}$.
\end{lem}

\medskip

\noi \textbf{Remarks. (i)} It is worth noting that theorem \ref{MainTheorem}  holds regardless of any causality assumption. Hypotheses $(1)$ and $(2)$ are pointwise energy conditions and condition $(3)$ a dynamic condition on the energy flow. This set of conditions is of a very different nature from the usual conditions of the classical singularity theorems.

\medskip

\textbf{(ii)} One cannot impose to $\wRic$ to be bounded, as Dajczer and Nomizu have proved in \cite{DajczerNomizu} that the spacetime is Einstein under this condition, so $\wRic$ is constant. The space $S^2_1\times\RR$, where $S^2_1$ is the $2$-dimensional Lorentz manifold of constant sectional curvature $1$, satisfies condition (2) and is not an Einstein manifold.

\medskip

\textbf{(iii)} The boundedness condition $(2)$ seems too demanding for theorem \ref{MainTheorem} to be of any physical interest. Our aim here is more to promote a method by a simple example rather than by a technical work which would improve these conditions. We shall nonetheless see in section \ref{Section2ndTheorem} how to weaken this assumption. The main probabilistic ingredient in the proof of theorem \ref{MainTheorem} is the basic explosion criterion stated in lemma \ref{LemmaExplosion}. This is a very general and rough tool which does not take into account any peculiar feature of our problem.

\medskip

\textbf{(iv)} Let $A$ be any subset of $\MMM$. The past domain of dependence of $A$ is the set $D^-(A)$ of points $m$ of $\MMM$ from which any future directed timelike path starting from $m$ eventually hits $A$. The future domain of dependence $D^+(A)$ is defined similarly using past directed timelike paths. The domain of dependence of $A$ is $D(A)= D^-(A)\cup D^+(A)$; it is globally hyperbolic. Note that if $A$ is a relatively compact spacelike hypersurface, and $m$ a point of $D^-(A)$, then there exists a constant $T(m)$, depending on $m$, \st all timelike path starting from $m$ hits $A$ before proper time $T(m)$.

Given $A\subset\MMM$, it is well-known that if the \textit{dominant energy condition} holds on its domain of dependence $D(A)$ and the energy momentum tensor vanishes on $A$ then it vanishes on the whole of $D(A)$, \cite{HawkingEllis}. Lemma \ref{LemmaGenerator} provides for free a similar result under different hypotheses. 

\begin{prop}
Let $A$ be a relatively compact spacelike hypersurface of a spacetime $(\MMM,g)$. Suppose the following conditions hold on $D^-(A)$:
\begin{itemize}
   \item strong and weak energy conditions,
   \item there exists a constant $c\in\RR$ \st $H_0\,\wRic\geq c\,\wRic$.
\end{itemize}
If $\Ric = 0$ on $A$ then $\Ric = 0$ on $D^-(A)$.
\end{prop}

Note that the strong and weak energy conditions together do not imply the dominant energy condition.

\ssk

\begin{Dem}
Combining the hypotheses of the proposition and lemma \ref{LemmaGenerator} we get the inequality
$$
\mcG\,\wRic \geq (\sigma^2+c)\,\wRic \geq 0,
$$
for a big enough $\sigma$. The function $\wRic$ is thus a non-negative $\mcG$-sub-harmonic function on $\OO\MMM$. Given $\Phi=(m,\be)$ with $m\in D^-(A)$, the stopping time $H=\inf\{s\geq 0\,;\,\pi\bigl(\Phi_s\bigr)\in A\}$ is \as bounded by a constant depending on $\Phi_0$, and 
$$
\EE_\Phi\bigl[\wRic(\Phi_H)\bigr] = \wRic(\Phi) + \EE_\Phi\Bigl[\int_0^H\bigl(\mcG\wRic\bigr)(\Phi_r)\,dr\Bigr] \geq \wRic(\Phi)\geq 0,
$$
by optional stopping. As $\wRic(\Phi_H)\equiv 0$ this proves that $\wRic(\Phi)=0$ and implies the result as $\Phi$ is any point of $\pi^{-1}\bigl(D^-(A)\bigr)$.
\end{Dem}

A similar result holds for $D^+(A)$ if $H_0\wRic\leq c\,\wRic$ for some constant $c\in\RR$.

\subsection{A second probabilistic singularity theorem}
\label{Section2ndTheorem}

We keep the notations of section \ref{SectionFirstTheorem}.

\medskip

\noi \textbf{a) Another explosion criterion.} We are going to use in this section a refined version of lemma \ref{LemmaExplosion} for which no boundedness hypothesis is needed.

\begin{lem}
\label{2ndExplosionLemma}
Suppose there exists two non-null, non-negative smooth functions on $\OO\MMM$, with $f\leq h$, and two constants $0\leq c'<c$ \st
$$
\mcG f\geq c\,f \quad \textrm{and} \quad \mcG h\leq c'\,h.
$$
Let $\Phi_0\in\OO\MMM$ be \st $f(\Phi_0)>0$. Then the basic relativistic diffusion started from $\Phi_0$ explodes with positive probability.
\end{lem}

\begin{Dem}
Suppose $\zeta$ is $\PP_{\Phi_0}$-\as infinite, so we can apply Ito's formula and write for any $s>0$
$$
f(\Phi_s) = f(\Phi_0) + M^f_s + \int_0^s\mcG f\,(\Phi_u)\,du,
$$
where $M^f$ is a $\PP_{\Phi_0}$-martingale \wrt the filtration generated by the diffusion process. Taking expectation and applying Fubini, we thus have 
$$
\EE_{\Phi_0}\bigl[f(\Phi_s)\bigr] = f(\Phi_0) + \int_0^s \EE_{\Phi_0}\bigl[\mcG f\,(\Phi_u)\bigr]\,du \geq f(\Phi_0) + c\int_0^s \EE_{\Phi_0}\bigl[f(\Phi_u)\bigr]\,du,
$$
so Gr\"onwall's lemma gives 
$$
\EE_{\Phi_0}\bigl[f(\Phi_s)\bigr] \geq f(\Phi_0)e^{c\,s}.
$$ 
Similarly, we have $\EE_{\Phi_0}\bigl[h(\Phi_s)\bigr] \leq~ h(\Phi_0)e^{c's}$. We get a contradiction noting that, since $f\leq h$, we should have 
$$
f(\Phi_0)e^{c\,s}\leq \EE_{\Phi_0}\bigl[f(\Phi_s)\bigr] \leq \EE_{\Phi_0}\bigl[h(\Phi_s)\bigr] \leq h(\Phi_0)e^{c's}
$$
for all times, which cannot happen.
\end{Dem}

We are going to apply this explosion criterion to the function $f=\wRic$ and a function $h$ of the form $\wRic + U$. We construct this function $U$ in the next paragraph.

\medskip

\noi \textbf{b) Green function on $\HH$.} Let denote by $\triangle$ the Laplacian on $\HH$ and by $G : \HH\times\HH\ra\RR^*_+$ the Green function of $\frac{1}{2}\triangle$. Due to the highly homogenous character of $\HH$, the quantity $G(x,y)$ is actually a function of the hyperbolic distance from $x$ to $y$. Let $x_{\textrm{ref}}\in\HH$ be any reference point and let $\rho$ denote the hyperbolic distance function to $x_{\textrm{ref}}$. A continuous function $f : \HH\ra \RR$ is said to be exp-bounded if $\sup_x\,e^{-a\rho^2(x)}\,\bigl|f(x)\bigr|<\infty$ for every $a>0$. This definition does not depend on the choice of reference point $x_{\textrm{ref}}$. It is well-known (consult for instance the classic book \cite{Friedman} of A. Friedman) that if $f$ is exp-bounded the equation 
$$
\frac{1}{2}\triangle u = -f
$$ 
has a unique solution which is null at infinity; it is given by the formula 
$$
u(x) = \int G(x,y)f(y)\,dy,
$$
where we write $dy$ for the Riemann volume form of $\HH$.

\ssk

Let us come back on $\OO\MMM$ and identify, for each $m\in\MMM$, the future unit tangent bundle $T^1_m\MMM$ to $\HH$ by arbitrarily identifying an element of $\OO_m\MMM$ to the canonical basis of $\RR^{1,3}$. The exp-boundedness character of the function $\wRic$ does not depend on this arbitrary choice; suppose it is exp-bounded for all $m\in\MMM$ and set, for $\Phi = \bigl(m,(\be_0,\cdots,\be_3)\bigr)$
\begin{equation}
\label{DefnU}
U(\Phi) = 2\int_{T^1_m\MMM}G(\be_0,y)\Ric_m(y,y)\,dy.
\end{equation}
As $G(\be_0,y)$ depends only on the hyperbolic distance from $\be_0$ to $y$ the function $U$ is well-defined independently of our arbitrary identifications. It solves the equation
\begin{equation}
\label{IdentityLaplacianU}
\frac{1}{2}\sum_{j=1}^3V_j^2\,U = -2\,\wRic,
\end{equation}
since the operator $\frac{1}{2}\sum_{j=1}^3\limits V_j^2$ on $\OO\MMM$ induces on each fiber $T^1_m\MMM$ the operator $\frac{1}{2}\triangle$.

\medskip

\noi \textbf{c) A second probabilistic singularity theorem.} The following singularity theorem is similar in nature to theorem \ref{MainTheorem}, and essentially states that a spacetime has a probabilistic singularity if some static and dynamical energy conditions hold. No causality assumption is needed. Set $h = \wRic + U$, and recall we write $R$ for the scalar curvature. 

\begin{thm}
\label{2ndMainTheorem}
Let $(\MMM,g)$ be a Lorentzian manifold satisfying the following conditions. 
\begin{enumerate}
   \item[(1')] \emph{Static energy conditions.} $\wRic$ is non-negative and non-identically null, and $R\leq 0$.
   \item[(2')] \emph{Regularity condition}. The function $\wRic$ is exp-bounded in each $T^1_m\MMM$ and there exists some constants $0<\al<1,\;0\leq c'<c$ and $\sqrt{\frac{c}{2}}<\sigma<\sqrt{\frac{c'}{2\,\al}}$, \st 
$$
\frac{1-\al}{\al}\,\wRic \leq U.
$$
   \item[(3')] \emph{Dynamical energy conditions.} \emph{(i)} $H_0\,\wRic \geq (c-2\,\sigma^2)\,\wRic$, 
   
$\quad\quad\quad\quad\quad\quad\quad\quad\quad\quad\quad\quad$\emph{(ii)} $H_0\,h \leq (c'-2\,\al\,\sigma^2)\,h$.
\end{enumerate}
Let $\Phi_0\in\OO\MMM$ be \st $\wRic(\Phi_0)>0$. Then the basic relativistic diffusion with diffusivity $\sigma^2$, started from $\Phi_0$, explodes with positive probability.
\end{thm}

\noi Note that our choice of constants gives $c'-2\,\al\,\sigma^2>0$ and $c-2\,\sigma^2<0$, which gives a non-trivial character to conditions (3').

\section{Conclusion}

We propose in this work a simple probabilistic method to probe certain aspects of the singular features of a spacetime, under the form of the existence of incomplete random dynamics obtained as random perturbation of the geodesic flow. Theorems \ref{MainTheorem} and \ref{2ndMainTheorem} show that geodesic completeness can be quite sensitive to random noise under some circumstances. This simple approach also highlights some unusual features, compared to the traditional studies on geodesic incompleteness.
\begin{itemize}
   \item[-] It might be interesting to work on the bundles $T^1\MMM$ or $\OO\MMM$ and not only on the base manifold $\MMM$.
   \item[-] The naive arguments given above indicate that causality and initial/boundary conditions might not be the crux of everything. In that respect, the appearance of a dynamic condition on the energy flow happens to be interesting, and needs to be compared to the usual pointwise/static energy conditions. 
\end{itemize}

\section{Proofs}
\label{SectionProof}

\subsection{Proofs of lemma \ref{LemmaLifting} and proposition \ref{PropTube}}

\subsubsection{Proof of lemma \ref{LemmaLifting}}

\noi $\bullet$ \textit{Existence.} Lift first arbitrarily the $\mcC^2$ path $\bigl(\gamma(s),\dot \gamma(s)\bigr)$ in $T^1\MMM$ into a $\mcC^1$ path $\Gamma_s = \bigl(\gamma(s),\be(s)\bigr)$ in $\OO\MMM$. As $\frac{d}{ds}\gamma(s)=\be_0(s)$, there exists some $\mcC^1$ real-valued controls $h^1,h^2,h^3$ and $\ell^{1,2},\ell^{1,3},\ell^{2,3}$ defined on $[0,T]$ \st 
$$
\dot \Gamma_s = H_0\bigl(\Gamma_s\bigr) + \sum_{j=1}^3\limits V_j\bigl(\Gamma_s\bigr)\,h^j(s) + \sum_{1\leq a<b\leq 3}\limits V_{ab}\bigl(\Gamma_s\bigr)\,\ell^{ab}(s),
$$ 
where $V_{ab}$ is the canonical vertical vector field on $\OO\MMM$ generated by the Lie element $E_{ab} = \ep_a\otimes\ep_b^*-\ep_b\otimes\ep_a^*$ of $SO(3)\subset SO_0(1,3)$. Let $A$ be the $SO_0(1,3)$-valued solution of the differential equation $dA_s =  -\sum_{1\leq a<b\leq 3}\limits A_sE_{ab}\,\ell^{ab}(s)$. Then the path $(\Psi_s)_{0\leq s\leq T} = (\Gamma_sA_s)_{0\leq s\leq T}$ satisfies the conditions of the lemma. 

\ssk

\noi $\bullet$ \textit{Uniqueness.} Suppoe $(\Theta_s)_{0\leq s\leq T}$ is another lift of $\bigl(\gamma(s),\dot\gamma(s)\bigr)_{0\leq s\leq T}$ to $\OO\MMM$ satisfying the above conditions with some controls $g^i$. As $\pi_1(\Theta_s)=\pi_1(\Psi_s)$, we need to have $\Theta_s = \Psi_s B_s$, for some $S0(3)\bigl(\subset S0_0(1,3)\bigr)$-valued $\mcC^2$ process $(B_s)_{0\leq s\leq T}$. Write $\dot B_s = \sum_{1\leq a<b\leq 3}\al^{ab}(s,B_s)V_{ab}(B_s)$, identifying here the vector fields $V_{ab}$ to vector fields on $SO(3)$. Then, we have on the one hand 
$$
\dot\Theta_s = H_0(\Theta_s) + \sum_{j=1}^3 V_j(\Theta_s)g^j(s),
$$
and on the other hand
$$
\dot\Theta_s = \dot\Psi_sB_s + \Psi_s\dot B_s = H_0(\Theta_s) + \sum_{j=1}^3 V_j(\Theta_s)h^j(s) + \sum_{1\leq a<b\leq 3}\al^{ab}(s,B_s)V_{ab}(\Theta_s)
$$
It follows that $g^j(s)=h^j(s)$, for each $j\in\{1,2,3\}$, and all the $\al^{ab}$ are identically null, so $\Theta_s=\Psi_sB_0$, for some $B_0\in SO(3)$, and eventually $\Theta_s=\Psi_s$, since $\Theta_0=\Psi_0=\Gamma_0$ .

\subsubsection{Proof of proposition \ref{PropTube}}

Using the chart $F$ on $\mcV$ provides a trivialization of $\OO\mcV$ for a small enough choice of $U$. We can thus consider $\OO\mcV$ as a submanifold of some open set of some $\RR^p$. Suppose $\Phi$ is in $\OO\mcV$. Then lemma \ref{LemmaLifting} and Stroock and Varadhan support theorem prove that 
$$
\PP_{\Phi}\Bigl((\Psi_s)_{0\leq s\leq T} \textrm{ exits }\mcV \textrm{ in }F\bigl(\{T\}\times U\bigr)\Bigr) > 0.
$$ 
To conclude in the general case it suffices to note that one can associate to any pair of points $(m_0,\dot m_0), (m_1,\dot m_1)$ in $T^1\MMM$ a timelike path $\rho : [0,1] \ra \MMM$ \st $\bigl(\rho(0),\dot\rho(0)\bigr) = (m_0,\dot m_0)$ and $\bigl(\rho(1),\dot\rho(1)\bigr) = (m_1,\dot m_1)$.

\subsection{Proofs of lemma \ref{LemmaGenerator} and theorem \ref{MainTheorem}}

\subsubsection{Proof of lemma \ref{LemmaGenerator}}

Write $\bigl(m,(\be_0,\be-1,\be_2,\be_3)\bigr)$ for a generic point of $\OO\MMM$. As each vector field $V_j$ induces no dynamics on $\MMM$ and generates in the $SO_0(1,3)$-fiber a hyperbolic rotation in the $2$-plane spanned by $\be_0$ and $\be_j$, we have 
\begin{equation*}
\begin{split}
\wRic(e^{tV_j}\Phi) &= - g\Bigl(R\bigl(\be_0,(\cosh t)\be_0+(\sinh t)\be_j\bigr)\bigl((\cosh t)\be_0+(\sinh t)\be_j\bigr),\be_0\Bigr)  \\
&+ \sum_{k=1}^3 g\Bigl(R\bigl(\be_k,(\cosh t)\be_0+(\sinh t)\be_j\bigr)\bigl((\cosh t)\be_0+(\sinh t)\be_j\bigr),\be_k\Bigr) \\
\end{split}
\end{equation*}
So we have
\begin{equation*}
V_j^2\wRic\,(\Phi) = \frac{d^2}{dt^2}_{\big|t=0}\,\wRic(e^{tV_j}\Phi) = 2\,\Ric_m(\be_0,\be_0) + 2\,\Ric_m(\be_j,\be_j),
\end{equation*}
and 
\begin{equation}
\label{IdentityLaplacianwRic} 
\sum_{j=1}^3 V_j^2\,\wRic = 4\,\wRic + 4\,\widetilde{{\bf T}},
\end{equation}
as $\sum_{j=1}^3\limits \Ric_m(\be_j,\be_j) = 2\,\widetilde{{\bf T}} - \Ric_m(\be_0,\be_0)$. The statement of the lemma follows.

\subsubsection{Proof of theorem \ref{MainTheorem}}
 
Let $\Phi_0$ be the starting point of the relativistic diffusion in $\OO\MMM$; we suppose $\wRic(\Phi_0)>0$. Under hypothesis $(2)$ we can combine lemmas \ref{LemmaExplosion} and \ref{LemmaGenerator} and see that the diffusion starting from $\Phi_0$ explodes with positive probability if there exists a positive constant $C$ \st $H_0\,\wRic + 2\,\sigma^2\,\wRic + 2\,\sigma^2\,\widetilde{{\bf T}} \geq C\,\wRic$. As $\widetilde{{\bf T}}\geq 0$ by hypothesis $(1)$, this will be the case if $H_0\,\wRic \geq \bigl(C-2\,\sigma^2\bigr)\,\wRic$; this condition is condition $(3)$.

\subsection{Proof of theorem \ref{2ndMainTheorem}}

The proof consists in checking that we can apply the explosion lemma \ref{2ndExplosionLemma} to the non-negative functions $f=\wRic$ and $h = \wRic + U\geq f$. As seen in lemma \ref{LemmaGenerator}, the condition $\mcG\,f\geq c\,f$ is equivalent to the inequality
\begin{equation}
\label{Cond1}
H_0\,\wRic + 2\,\sigma^2\,\wbfT\geq (c-2\,\sigma^2)\,\wRic;
\end{equation}
it follows from  condition (3'-i) since $\wbfT\geq 0$. The condition $\mcG h\leq c'h$ reads 
$$
H_0\,h + \frac{\sigma^2}{2}\sum_{j=1}^3V_j^2\wRic +2\,\sigma^2\,\wbfT + \frac{\sigma^2}{2}\sum_{j=1}^3V_j^2 U \leq c'h.
$$
By \eqref{IdentityLaplacianwRic} and \eqref{IdentityLaplacianU}, it is equivalent to 
\begin{equation}
\label{Cond2}
H_0 h + 2\,\sigma^2\,\wbfT \leq  c'\,h.
\end{equation}
To see that \eqref{Cond2} follows from condition (3'-ii) notice that the inequality $\wbfT\leq \al\,h$ is equivalent to the inequality $\frac{1-\al}{\al}\,\wRic + \frac{R}{2\,\al}\leq U$. This condition is implied by condition (2') as we suppose $R\leq 0$.


\def\cprime{$'$} \def\cprime{$'$}

\end{document}